\documentclass{article}
\usepackage{amssymb}
\usepackage{graphicx}
\usepackage{amsmath}

\newtheorem{theorem}{Theorem}
\newtheorem{acknowledgement}[theorem]{Acknowledgement}
\newtheorem{algorithm}[theorem]{Algorithm}

\newtheorem{corollary}[theorem]{Corollary}

\newtheorem{example}[theorem]{Example}

\newtheorem{lemma}[theorem]{Lemma}

\newtheorem{remark}[theorem]{Remark}

\newenvironment{proof}[1][Proof]{\textbf{#1.} }{\ \rule{0.5em}{0.5em}}
\input{tcilatex}
\numberwithin{theorem}{section}

\begin{document}

\title{Algorithms for Producing and Ordering Lexical and Nonlexical
Sequences out of one Element\\
\bigskip }
\author{Elias Abboud\bigskip \\
{\small Beit Berl College, Doar Beit Berl, 44905 Israel }{\footnotesize 
\thanks{%
Also the Academic Arab College - Haifa.}}\\
{\small \ Email: eabboud @beitberl.ac.il}}
\maketitle

\begin{abstract}
This paper deals with algorithms for producing and ordering lexical and
nonlexical sequences of a given degree. The notion of \textquotedblright
elementary operations\textquotedblright\ on positive $\alpha $-sequences is
introduced. Our main theorem answers the question of when two lexical
sequences are adjacent. Given any lexical sequence, $\alpha \in L_{n},$ we
can produce its adjacent successor as follows; apply one elementary
operation on the tail of the longest left sequence, of even length, which
gives a lexical successor $\alpha ^{\prime }\in L_{n}$, then compute the
fundamental sequence $f=\alpha \wedge \alpha ^{\prime }$ $\in L_{m}$ and
conclude for $m\nmid n$ $\ $that $\alpha $ is adjacent to $\alpha ^{\prime }$
in $L_{n}.$ Whereas for $m\mid n$ , the sequence $\alpha $ is adjacent to a
sequence generated by $\ f$ and the least element of $L_{d},$ where $d=\frac{%
n}{m}.$ Thus, while right sequences control the lexicality property of an $%
\alpha $-sequence, it turns out that left sequences control the adjacency
property of lexical and nonlexical sequences.

\bigskip

 \smallskip lexical sequences, adjacent.

\bigskip

\end{abstract}

\section{Introduction}

The topic of lexical sequences is originated from the field of Discrete
Dynamical Systems, and their combinatorial properties are related to
patterns of words on two letters which label functions constituting the
inverse graph $G_{\zeta }^{(n)}$ of the $n$ th iterate of the parabolic map $%
p_{\zeta }(x)=\zeta x(2-x),$ $\zeta \in \lbrack 0,\infty ),$ $x\in \lbrack
0,\infty ),($ see \cite{BLMS}, \cite{L}, \cite{LS}).

A lexical order denoted by ''$<$'' is a total order defined on the set of
positive $\alpha $-sequences, which was introduced in \cite{LM} and was
extended in \cite{BLMS} and used by others in studying maps of the interval,
(see \cite{B}, \cite{HZ}, \cite{SH}).

This paper deals with algorithms for producing and ordering lexical and
nonlexical sequences of a given degree. A problem about the parity of
adjacent lexical sequences was posed in a paper by J. D. Louck \cite{L}.
This problem was solved in \cite{CLW}, and the authors proved that the
length of adjacent lexical sequences have opposite parity. They used a
transitivity property, as a strategy, in order to avoid the direct question
of when two sequences are adjacent as they clarified in their paper, \cite[%
p.484]{CLW}. Their method is resumed as follows; given two adjacent lexical
sequences they proved a ''structural theorem'', expressed them by means of
their ''fundamental'' sequence and then inferred the opposite parity
property.

On the contrary, in this paper we do attack this direct question. Given one
lexical sequence, we will be able to list all its successors and all its
precedings.

The paper is organized as follows; after recalling the common terminology
needed for this work, we introduce in section \ref{EO}, the notion of
\textquotedblright elementary operations\textquotedblright\ on positive $%
\alpha $-sequences and their properties. These will be applied on the tail
of the longest left sequence in the set of all $\alpha $-sequences and on
the tail of the longest left sequence which give a lexical successor in the
set of all lexical sequences. In section \ref{ACLS}, we will prove the main
theorem which answers the question of when two lexical sequences are
adjacent. Precisely, given any lexical sequence, $\alpha \in L_{n},$ we can
produce its adjacent successor as follows; apply one elementary operation on
the tail of the longest left sequence, of even length, which gives a lexical
successor $\alpha ^{\prime }\in L_{n}$, compute the fundamental sequence $%
f=\alpha \wedge \alpha ^{\prime }$ $\in L_{m}$ and conclude for $m\nmid n$ $%
\ $that $\alpha $ is adjacent to $\alpha ^{\prime }$ in $L_{n},$ while for $%
m\mid n$ , the sequence $\alpha $ is adjacent to a sequence generated by $\
f $ and the least element of $L_{d},$ where $d=\frac{n}{m}.$

Finally, in section \ref{Algthms}, we will produce and order the main sets
of lexical and nonlexical sequences out of one element.

\subsection{\protect\bigskip Definitions and Notations}

We borrow some of the terminology that appeared in the references \cite{BLMS}%
, \cite{CLW}, \cite{L} and \cite{LM}. The principal notations will be
adopted from \cite{CLW}.

\begin{itemize}
\item An\textit{\ }$\alpha $\textit{-sequence} is a finite sequence over the
natural numbers $\mathbb{N=}\left\{ 1,2,3,...\right\} $; $\alpha =(\alpha
_{1},\alpha _{2},...,\alpha _{k})$ where, $\alpha _{i}\in \mathbb{N} $. A 
\textit{right sequence} of $\alpha =(\alpha _{1},\alpha _{2},...,\alpha
_{k}) $ is any sequence $(\alpha _{i},\alpha _{i+1},...,\alpha _{k}),$ $%
1\leq i\leq k$. \ A \textit{left sequence} of \ $\alpha $ is any sequence of
the form $(\alpha _{1},\alpha _{2},...,\alpha _{i}),$ $1\leq i\leq k.$ In
this last case $\alpha _{i}$ is called the \textit{tail} of the left
sequence.

\item A \textit{total order} relation on the set of all \ $\alpha $%
-sequences is defined as follows: With each $\alpha $-sequence $\tau =(\tau
_{1},\tau _{2},$...,$\tau _{k})$ we associate the sequence 
\begin{equation*}
Al(\tau )=(\tau _{1},-\tau _{2},...,(-1)^{k}\tau _{k},0,0,...)
\end{equation*}
\end{itemize}

which alternates in sign of its nonzero elements. For any different $\alpha $%
-sequences $\tau $ and $\tau ^{\prime }$ we form the difference 
\begin{equation*}
Al(\tau )-Al(\tau ^{\prime }).
\end{equation*}

If the first nonzero term in this difference is positive, we write $\tau
>\tau ^{\prime }$ and if, otherwise, is negative we write $\tau <\tau
^{\prime }.$

\begin{itemize}
\item If $\alpha <\beta $ then we call $\beta $ a \textit{successor} of $%
\alpha ,$ and $\alpha $ a \textit{preceding} of $\beta .$

\item Two sequences $\alpha $ and $\beta $ in the same set $S$ are \textit{%
adjacent} in $S$ if there is no $\gamma \in S$ such that $\alpha <\gamma
<\beta .$

\item We say that $\alpha $ is \textit{adjacent to }$\beta $ in $S,$ and
denote $\alpha \underset{adj}{<}\beta ,$ if $\alpha $ and $\beta $ are 
\textit{adjacent} in $S$ and $\alpha <\beta .$

\item A sequence $\alpha =(\alpha _{1},\alpha _{2},...,\alpha _{k})$ is
called \textit{lexical \ }iff\textit{\ }$\alpha >(\alpha _{i},\alpha
_{i+1},...,\alpha _{k}),$ for all right sequences, $2\leq i\leq k$,
otherwise is called \textit{nonlexical}.

\item \ $l(\alpha )$ denotes the \textit{length} of the $\alpha $ and is
defined as the number of its elements $l(\alpha )=k.$

\item The \textit{degree} $D(\alpha )$ of the sequence $\alpha $ is defined
by $D(\alpha )=\sum_{i=1}^{k}\alpha _{i}.$

\item The set of positive sequences of a given degree $\ \ \ \ \ \ \ \ \ $%
\begin{equation*}
\ A_{n}=\left\{ (\alpha _{1},\alpha _{2},...,\alpha _{k})\right| \alpha
_{i}\in \mathbb{N},\sum_{i=1}^{k}\alpha _{i}=n,\ 1\leq k\leq n\}.
\end{equation*}

\item The set of lexical sequences of a given degree $\ \ \ \ \ \ \ \ \ \ \
\ \ \ \ \ \ \ \ \ \ \ \ \ \ \ \ \ \ \ \ \ \ \ \ \ \ \ \ \ \ \ \ $%
\begin{equation*}
L_{n}=\left\{ \alpha \right| \alpha \text{ is lexical},1+D(\alpha )=n\}.
\end{equation*}

\item The set of lexical sequences of degree dividing $n.$ 
\begin{equation*}
D_{n}=\underset{d\mid n}{\cup }L_{d}
\end{equation*}

\item The \textit{product} of two sequences $\alpha \beta ,$ where $\alpha
=(\alpha _{1},\alpha _{2},...,\alpha _{k})$ and $\beta =(\beta _{1},\beta
_{2},...,\beta _{k})$ is defined as $\alpha \beta =(\alpha _{1},\alpha
_{2},...,\alpha _{k},\beta _{1},\beta _{2},...,\beta _{k}).$ In the same
manner, $\alpha ^{q}$ denotes the product of $\alpha $ with itself $q$ times
and $1^{q}=(1)^{q}$ i.e. $1$ is repeated $q$ times.

\item The sequence $\alpha \wedge \beta $ constitutes the longest common
left factor of $\alpha $ and $\beta .$ More precisely, $\alpha \wedge \alpha
=\alpha $ and for $\alpha \neq $ $\beta $ if $i$ is the first place for
which $\alpha _{i}\neq $ $\beta _{i}\,$then, $\alpha \wedge \beta =(\alpha
_{1},\alpha _{2},...,\alpha _{i-1},\alpha _{i})$ for $\alpha _{i}<$ $\beta
_{i}$ and $\alpha \wedge \beta =(\alpha _{1},\alpha _{2},...,\alpha
_{i-1},\beta _{i})$ for $\beta _{i}<\alpha _{i}.$

\item The sequences $\alpha _{e}$ of even length and $\alpha _{o}$ of odd
length are defined as follows:
\end{itemize}

\begin{equation*}
\alpha _{e}=\left\{ 
\begin{array}{c}
(\alpha _{1},\alpha _{2},...,\alpha _{k}+1),\text{ \ }l(\alpha )\text{ even}
\\ 
(\alpha _{1},\alpha _{2},...,\alpha _{k},1),\text{ \ }l(\alpha )\text{ odd}%
\end{array}
\right.
\end{equation*}
\begin{equation*}
\ \alpha _{o}=\left\{ 
\begin{array}{c}
(\alpha _{1},\alpha _{2},...,\alpha _{k},1),\text{ \ }l(\alpha )\text{ even}
\\ 
(\alpha _{1},\alpha _{2},...,\alpha _{k}+1),\text{ \ }l(\alpha )\text{ odd}%
\end{array}
\right.
\end{equation*}

\begin{itemize}
\item The zero and the first \textit{harmonic} of a sequence $\alpha $ are
defined as $h_{0}(\alpha )=\alpha ,h_{1}(\alpha )=\alpha _{o}\alpha .$ The $%
j $-th harmonic of $\alpha $ is defined recursively by 
\begin{equation*}
h_{j}(\alpha )=(h_{j-1}(\alpha ))_{o}h_{j-1}(\alpha )\text{ , }j=1,2,...
\end{equation*}

\item A sequence is \textit{fundamental} if it cannot be written as the
harmonic of another sequence.

\item The \textit{star product}, $\alpha \star \beta $, of two positive
sequences $\alpha $ and $\beta =(\beta _{1},\beta _{2},...,\beta _{k})$ is
defined as follows 
\begin{equation*}
\alpha \star \beta =\alpha _{o}(\alpha _{e})^{\beta _{1}-1}\alpha
_{o}(\alpha _{e})^{\beta _{2}-1}...\alpha _{o}(\alpha _{e})^{\beta
_{k}-1}\alpha
\end{equation*}
\end{itemize}

\section{Elementary Operations\label{EO}}

Let $\alpha =(\alpha _{1},\alpha _{2},...,\alpha _{k})$. We consider $\alpha 
$ as a sequence of $k$ \textit{cells }$C_{1}(\alpha ),...,C_{k}(\alpha )$
with \textit{values}: $val(C_{j}(\alpha ))=\alpha _{j}.$ We also consider $%
C_{1}(\alpha ),C_{3}(\alpha ),...$ as \textit{positive }cells, whereas $%
C_{2}(\alpha ),C_{4}(\alpha ),...$ as \textit{negative }cells.

Thus, $\alpha _{1},\alpha _{3},...$ are placed in positive cells whereas $%
\alpha _{2},\alpha _{4},...$ are placed in\textit{\ }negative \textit{cells. 
}We define the following two elementary operations on $\alpha -$sequences.

\begin{itemize}
\item \textit{Splitting cells.}
\end{itemize}

We split a cell into two adjacent cells, the first with value $\alpha _{i}-1$
and the second with value $1$. Formally, \textit{for \ }$\alpha _{i}>1$%
\textit{\ define, }

\begin{equation*}
\mathcal{E}_{s}:\text{ \ \ \ \ \ \ \ \ }\alpha ^{\prime }=\left\{ 
\begin{array}{c}
(\alpha _{1},...,\alpha _{i}-1)_{o}(\alpha _{i+1},...,\alpha _{k}).\text{ \
\ \ \ \ \ \ \ }l(\alpha _{1},...,\alpha _{i})\mathit{\ \ }\text{\textit{even}%
\ \ } \\ 
(\alpha _{1},...,\alpha _{i}-1)_{e}(\alpha _{i+1},...,\alpha _{k}).\text{ \
\ \ \ \ \ \ \ }l(\alpha _{1},...,\alpha _{i})\mathit{\ \ }\text{\textit{odd}}%
\end{array}
\right.
\end{equation*}

\begin{itemize}
\item \textit{Conjugation of cells.}
\end{itemize}

We conjugate a cell with value $\alpha _{i}=1$ and its left adjacent cell
with value $\alpha _{i-1}$ into one cell with value $\alpha _{i-1}+1.$
Formally, \textit{if }$\alpha _{i}=1$\textit{\ define, } 
\begin{equation*}
\mathcal{E}_{c}:\text{ \ \ \ \ \ \ \ \ }\alpha ^{\prime }=\left\{ 
\begin{array}{c}
(\alpha _{1},...,\alpha _{i-1})_{o}(\alpha _{i+1},...,\alpha _{k})\mathit{.}%
\text{ \ \ \ \ \ \ \ \ \ \ }l(\alpha _{1},...,\alpha _{i})\mathit{\ \ }\text{%
\textit{even}} \\ 
(\alpha _{1},...,\alpha _{i-1})_{e}(\alpha _{i+1},...,\alpha _{k})\mathit{.}%
\text{ \ \ \ \ \ \ \ \ \ \ \ \ }l(\alpha _{1},...,\alpha _{i})\mathit{\ \ }%
\text{\textit{odd}}%
\end{array}
\right.
\end{equation*}

\begin{example}
\label{example L7}The following example gives a motivation for splitting and
conjugation of cells. By starting from the minimal element in $L_{7}$,
applying $\mathcal{E}_{s}$ and $\mathcal{E}_{c}$ on negative cells only we
may produce all lexical sequences in $L_{7}$ successively from down and up.
\end{example}

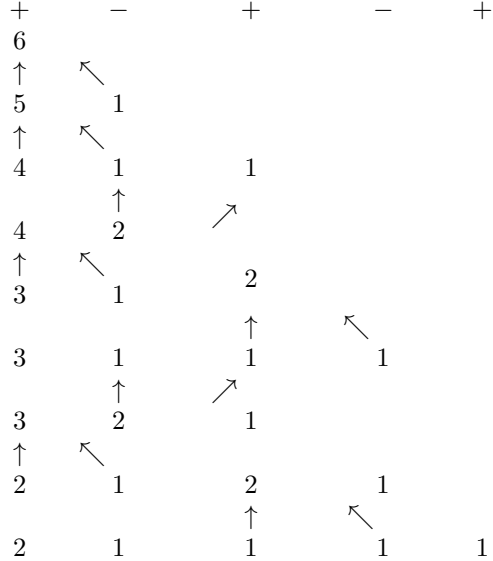
\begin{equation*}
\begin{array}{ccccc}
+ & - & + & - & + \\ 
6 &  &  &  &  \\ 
\begin{array}{c}
\uparrow \\ 
5%
\end{array}
& 
\begin{array}{c}
\nwarrow \qquad \\ 
1%
\end{array}
&  &  &  \\ 
\begin{array}{c}
\uparrow \\ 
4%
\end{array}
& 
\begin{array}{c}
\nwarrow \qquad \hspace{0in} \\ 
1%
\end{array}
& 
\begin{array}{c}
\\ 
1%
\end{array}
&  &  \\ 
\begin{array}{c}
\\ 
4%
\end{array}
& 
\begin{array}{c}
\uparrow \\ 
2%
\end{array}
& 
\begin{array}{c}
\nearrow \quad \ \ \ 
\end{array}
&  &  \\ 
\begin{array}{c}
\uparrow \\ 
3%
\end{array}
& 
\begin{array}{c}
\nwarrow \qquad \\ 
1%
\end{array}
& 2 &  &  \\ 
\begin{array}{c}
\\ 
3%
\end{array}
& 
\begin{array}{c}
\\ 
1%
\end{array}
& 
\begin{array}{c}
\uparrow \\ 
1%
\end{array}
& 
\begin{array}{c}
\nwarrow \qquad \\ 
1%
\end{array}
&  \\ 
\begin{array}{c}
\\ 
3%
\end{array}
& 
\begin{array}{c}
\uparrow \\ 
2%
\end{array}
& 
\begin{array}{c}
\nearrow \qquad \\ 
1%
\end{array}
& 
\begin{array}{c}
\end{array}
& 
\begin{array}{c}
\end{array}
\\ 
\begin{array}{c}
\uparrow \\ 
2%
\end{array}
& 
\begin{array}{c}
\nwarrow \qquad \\ 
1%
\end{array}
& 
\begin{array}{c}
\\ 
2%
\end{array}
& 
\begin{array}{c}
\\ 
1%
\end{array}
& 
\begin{array}{c}
\end{array}
\\ 
\begin{array}{c}
\\ 
2%
\end{array}
& 
\begin{array}{c}
\\ 
1%
\end{array}
& 
\begin{array}{c}
\uparrow \\ 
1%
\end{array}
& 
\begin{array}{c}
\nwarrow \ \quad \  \\ 
1%
\end{array}
& 
\begin{array}{c}
\\ 
1%
\end{array}%
\end{array}%
\end{equation*}

\begin{remark}
Starting from the maximal element $(6)\in L_{7}$, we may apply $\mathcal{E}%
_{s}$ and $\mathcal{E}_{c}$ on positive cells only and produce all lexical
sequences in $L_{7}$ successively from up and down.
\end{remark}

\subsection{Basic Properties of the Elementary Operations \label{Bproperties}%
}

It turns out that the elementary operations play an important role in
answering the question of ''when two $\alpha -sequences$ are adjacent?''

First, we answer the question for the set $A_{n}.$ But we emphasize the fact
that this work was inspired by the opposite parity property of adjacent
sequences which was posed in \cite{L} and solved in \cite{CLW}. The next
trivial lemma gives a motivation.

\begin{lemma}
Suppose that $\alpha ^{\prime }\in A_{n}$ is obtained from $\alpha \in A_{n}$
by applying one elementary operation $\mathcal{E}_{s}$ or $\mathcal{E}_{c}$
on a negative or positive cell of $\alpha ,$ then $\alpha $ and $\alpha
^{\prime }$ have opposite parity of their length.
\end{lemma}

Next, we show that applying an elementary operation gives a successor or a
preceding element.

\begin{lemma}
\begin{enumerate}
\item Suppose that $\alpha ^{\prime }\in A_{n}$ is obtained from $\alpha \in
A_{n}$ by applying one elementary operation $\mathcal{E}_{s}$ or $\mathcal{E}%
_{c}$ on a negative cell of $\alpha ,$ then $\alpha <$ $\alpha ^{\prime }.$

\item If $\alpha ^{\prime }\in A_{n}$ is obtained from $\alpha \in A_{n}$ by
applying one elementary operation $\mathcal{E}_{s}$ or $\mathcal{E}_{c}$ on
a positive cell of $\alpha ,$ then $\alpha ^{\prime }<$ $\alpha .$
\end{enumerate}
\end{lemma}

\begin{proof}
The elementary operations $\mathcal{E}_{s}$, $\mathcal{E}_{c}$ increase or
decrease the value of the corresponding cell according whether it is a
positive or a negative cell.
\end{proof}

In the following lemma we shall apply one elementary operation $\mathcal{E}%
_{s}$ or $\mathcal{E}_{c}$ on the tail of the longest left sequence of even
length, which is the first negative cell of $\alpha $ from the right. Also,
we shall apply one elementary operation $\mathcal{E}_{s}$ or $\mathcal{E}%
_{c} $ on the tail of the longest left sequence of odd length, which is the
first positive cell of $\alpha $ from the right.

\begin{lemma}
\label{adj-An}

\begin{enumerate}
\item Suppose that $\alpha ^{\prime }\in A_{n}$ is obtained from $\alpha \in
A_{n}$ by applying one elementary operation $\mathcal{E}_{s}$ or $\mathcal{E}%
_{c}$ on the tail of the longest left sequence, of $\alpha ,$ of even
length. Then, $\alpha \underset{adj}{<}\alpha ^{\prime }.$

\item Suppose that $\alpha ^{\prime }\in A_{n}$ is obtained from $\alpha \in
A_{n}$ by applying one elementary operation $\mathcal{E}_{s}$ or $\mathcal{E}%
_{c}$ on the tail of the longest left sequence, of $\alpha ,$of odd length.
Then, $\alpha ^{\prime }\underset{adj}{<}\alpha .$
\end{enumerate}
\end{lemma}

\begin{proof}
We shall prove only the first part for $\alpha _{i}>1$, the rest will follow
similarly.

Suppose that $\alpha =(\alpha _{1},\alpha _{2},...,\alpha _{i-1},\alpha
_{i},\alpha _{i+1}),$ where $i$ is even and $\alpha _{i+1}\geq 0.$ The
longest left sequence of $\alpha ,$ of even length, is $(\alpha _{1},\alpha
_{2},...\alpha _{i}).$ Thus, we have $\alpha ^{\prime }=(\alpha _{1},\alpha
_{2},...\alpha _{i}-1,1,\alpha _{i+1}).$ Let $\beta \in A_{n}$ satisfies, $%
\alpha \leq \beta \leq \alpha ^{\prime }.$

Case I: $\beta =(\alpha _{1},\alpha _{2},...\alpha _{i}-1,1)\gamma ,$ such
that $\gamma =(\gamma _{1},...,\gamma _{j})\in A_{\alpha _{i+1}}$. But
since, $\alpha ^{\prime }$ succeeds $\beta ,$ we must have $j=1$ and $\gamma
=\alpha _{i+1}.$ So, $\gamma =\alpha ^{\prime }.$

Case II: $\beta =(\alpha _{1},\alpha _{2},...\alpha _{i})\gamma ,$ such that 
$\gamma =(\gamma _{1},...,\gamma _{j})\in A_{\alpha _{i+1}}$. Here again
since, $\beta $ succeeds $\alpha ,$ we must have $j=1$ and $\gamma =\alpha
_{i+1}.$ So, $\gamma =\alpha .$ Consequently, $\alpha $ is adjacent to $%
\alpha ^{\prime }$ in $A_{n}.$
\end{proof}

The following lemma covers a subcase of Theorem \ref{adjcrth}, \ in section %
\ref{ACLS}.

\begin{lemma}
\label{gosigma}Let $g\in L_{m}$ be fundamental such that for some positive
sequences $\sigma $, $\kappa $ we have $\alpha =g_{e}\sigma $ and $\beta
=g_{o}\kappa $ are adjacent in $L_{2m}$, then $g_{o}\sigma \in $ $L_{2m}.$

\begin{proof}
By Theorem 3.4 (2) \cite[p. 1994]{DLW}, $\ $we have $\ \alpha =g_{e}\sigma
=g_{e}\widehat{g}$ and $\beta =g\star \lambda _{2}$, where $\lambda _{2}$ is
the least element of $L_{2}.$ and the structures of $g$ and $\widehat{g}$
are given by one of two possibilities (see the proof of Lemma 3.3 \cite[p.
1993]{DLW} and its preceding paragraph );

(i) There is a $\tau \in L_{m_{1}},$ where $m_{1}\mid m,$ and an odd integer 
$r>1$ such that $g=\tau \star \lambda _{r}=\tau \star (2,1^{r-3})=\tau
_{o}\tau _{e}\tau _{o}^{r-3}\tau $ and $\widehat{g}=\tau \star
(1^{r-1})=\tau _{o}^{r-1}\tau .$ Since $\tau _{o}\tau _{e}\tau _{o}^{r-3}$
has odd parity, we have $g_{o}=\tau _{o}\tau _{e}\tau _{o}^{r-3}\tau _{e}$.
Thus 
\begin{equation*}
g_{o}\sigma =\ g_{o}\widehat{g}=\tau _{o}\tau _{e}\tau _{o}^{r-3}\tau
_{e}\tau _{o}^{r-1}\tau =\left\{ 
\begin{array}{c}
(\tau _{o}\tau _{e})\tau _{o}^{r-4}(\tau _{o}\tau )\tau _{o}^{r-1}\tau ,%
\text{\ \ \ \ \ for odd }r\geq 5 \\ 
\tau _{o}\tau _{e}^{2}\tau _{o}^{2}\tau ,\text{ \ \ \ \ \ \ for }r=3\text{\
\ \ }%
\end{array}%
.\right.
\end{equation*}

Equivalently,%
\begin{equation*}
g_{o}\sigma =\left\{ 
\begin{array}{c}
\tau \star (2,1^{r-4},2,1^{r-1}),\text{ \ \ \ \ \ \ \ \ \ \ \ for odd }r\geq
5 \\ 
\tau \star (3,1,1),\text{ \ \ \ \ \ \ \ for }r=3\text{\ \ \ }%
\end{array}%
,\right.
\end{equation*}

where for an empty word $\tau $, 
\begin{equation*}
g_{o}\sigma =\left\{ 
\begin{array}{c}
(2,1^{r-4},2,1^{r-1}),\text{ \ \ \ \ \ \ \ \ \ \ \ for odd }r\geq 5 \\ 
(3,1,1),\text{ \ \ \ \ \ \ \ for }r=3\text{\ \ \ }%
\end{array}%
.\right.
\end{equation*}

Since $\tau ,(2,1^{r-4},2,1^{r-1}),$ $(3,1,1)$ are lexical and the star
product of lexical sequences is lexical ( Lemma 3.3 \cite[p. 491]{CLW}), the
result follows.

(ii) There is a $\tau \in L_{m_{1}},$ with $m_{1}\nmid m$ and a positive
sequence $\zeta $ such that $g=\tau _{o}\zeta $ and $\widehat{g}=\tau
_{e}\zeta ,$ where $g$ is not least in $L_{m}$. By Lemma 3.3 (2) \cite[p.1993%
]{DLW}, $\widehat{g}=\tau _{e}\zeta $ is adjacent to $g=\tau _{o}\zeta $ in $%
L_{m}.$ Whereas, by Lemma 2.5 \cite[p.1988]{DLW}, $g_{e}\widehat{g}$ $\in
L_{2m}.$ Following the proof of Lemma 2.5 \cite[p.1988]{DLW}, we claim that $%
g_{o}\sigma =\ g_{o}\widehat{g}$ is lexical. Otherwise, there exists a right
factor $\mu $ of $g$ such that $\mu _{p}\widehat{g}>g_{o}\widehat{g},$ where 
$\mu _{p}=\mu _{e}$ or $\mu _{o}.$ The existence of such $\mu $ is explained
by the fact that any right factor $\rho $ of $\widehat{g}$ satisfies $\rho <$
$\widehat{g}<$ $g_{e}\widehat{g}<$ $g_{o}\widehat{g}$. Hence, $\mu _{p}%
\widehat{g}>g_{o}\widehat{g}>\mu ,$ which yields that $\mu _{p}=$ $\mu _{o}$
is odd and $g=\mu _{o}g^{\prime }$ for some right factor $g^{\prime }$ of $g.
$ If $\mu _{o}$ is also a left factor of $\widehat{g}$, then $\widehat{g}%
=\mu _{o}\widehat{g}^{\prime }$ and by the lexicality of $\widehat{g}$ we
must have $\widehat{g}>\widehat{g}^{\prime }$. But since $\mu _{o}$ is of
odd parity, we get $\widehat{g}=\mu _{o}\widehat{g}^{\prime }>\mu _{o}%
\widehat{g}>g_{o}\widehat{g}>g$ which is a contradiction. If $\mu _{o}$ is
not a left factor of $\widehat{g}$ then $l(\mu )\geq l(\tau )$ which implies
that $g^{\prime }$ is a right factor of $\zeta $ and thus for $\widehat{g}.$
By the lexicality of $\widehat{g}$ we must have $\widehat{g}>g^{\prime }$
and so $g=\mu _{o}g^{\prime }>\mu _{o}\widehat{g}>g_{o}\widehat{g}>g$ which
is a contradiction. This proves indeed that $\ g_{o}\sigma =\ g_{o}\widehat{g%
}$ is lexical.
\end{proof}
\end{lemma}

\section{\protect\bigskip Adjacency Criterion for Lexical Sequences\label%
{ACLS}}

Now we are ready to deal with the lexical sequences in $L_{n}$ and $D_{n}.$
We shall use theorems and lemmas from \cite{CLW} for proving our main
theorem. We shall apply one elementary operation on the cell corresponding
to the tail of the longest left sequence of even length, which gives a
lexical successor, which is equivalently, the first negative cell of the
sequence from the right, which gives a lexical successor. Similarly, the
cell corresponding to the tail of the longest left sequence of odd length
which gives a lexical preceding, is the same as the first positive cell of $%
\alpha ,$ from the right, which gives a lexical preceding.

\begin{theorem}
\label{adjcrth}Let $\alpha \in L_{n}$ and suppose that $\alpha $ is not
maximal in $L_{n}.$ Suppose that $\alpha ^{\prime }\in L_{n}$ is obtained
from $\alpha $ by applying one elementary operation $\mathcal{E}_{s}$ or $%
\mathcal{E}_{c}$ on the cell corresponding to the tail of the longest left
sequence of even length, which gives a lexical successor in $L_{n}$. Let $f=$
$\alpha $ $\wedge \alpha ^{\prime }$. Assume further that, $f\in L_{m}$ and $%
n=md+r$, $0\leq r<m.$

\begin{enumerate}
\item If $\ r>0,$ then $\alpha $ \ is adjacent to $\alpha ^{\prime }$ in $%
L_{n}.$

\item If $\ r=0,$ then $\alpha $ \ is adjacent to $f\star \lambda $ in $%
L_{n},$ where $\lambda $ is the least element of $L_{d}$.
\end{enumerate}
\end{theorem}

\begin{example}
In $L_{11},$ the longest left sequence of even length, which gives a lexical
successor of $\alpha =(3,2,3,2)$ when applying $\mathcal{E}_{s}$ is $(3,2).$%
\begin{equation*}
\begin{array}{cccccc}
& + & - & + & - & + \\ 
\alpha ^{\prime } & 3 & 1 & 1 & 3 & 2 \\ 
\begin{array}{c}
\\ 
\alpha%
\end{array}
& 
\begin{array}{c}
\\ 
3%
\end{array}
& 
\begin{array}{c}
\uparrow \\ 
2%
\end{array}
& 
\begin{array}{c}
\nearrow \quad \quad \\ 
3%
\end{array}
& 
\begin{array}{c}
\\ 
2%
\end{array}
& 
\end{array}
.
\end{equation*}
In this case $\alpha $ is adjacent to $\alpha ^{\prime }$ in $L_{11}.$
\end{example}

\begin{proof}
(of Theorem \ref{adjcrth})

Suppose that the cell corresponding to the tail of the longest left sequence
of even length, which gives a lexical successor, has value $a.$ Denote this
left sequence of even length by $\gamma (a).$ Obviously, $\gamma $ has odd
length. For $a>1$ we are applying the elementary operation $\mathcal{E}_{s}$
and for $a=1$ we are applying $\mathcal{E}_{c}.$ Let 
\begin{equation}
\alpha =\gamma (a)\rho \in L_{n}  \label{1}
\end{equation}
then, 
\begin{equation}
\alpha ^{\prime }=\left\{ 
\begin{array}{c}
\gamma (a-1,1)\rho ,\text{ \ \ for \ \ }a>1 \\ 
\gamma _{o}\rho ,\text{ \ \ \ \ \ \ \ \ \ \ \ \ \ \ for \ }a=1%
\end{array}
\right. \text{, \ \ }\alpha ^{\prime }\in L_{n}.  \label{2}
\end{equation}
Clearly we have, 
\begin{equation}
f=\alpha \wedge \alpha ^{\prime }=\left\{ 
\begin{array}{c}
\gamma (a-1)\text{ \ \ for \ \ }a>1 \\ 
\gamma \text{\ \ \ \ \ \ \ \ \ \ \ \ \ \ \ for \ }a=1%
\end{array}
\right. \text{, \ }f\in L_{m}.  \label{3}
\end{equation}
Therefore, 
\begin{equation*}
f_{e}=\left\{ 
\begin{array}{c}
\gamma (a)\text{ \ \ \ \ \ \ \ \ \ \ \ \ for \ \ }a>1 \\ 
\gamma (1)\text{\ \ \ \ \ \ \ \ \ \ \ \ \ \ for \ }a=1%
\end{array}
\right.
\end{equation*}
equivalently, 
\begin{equation}
f_{e}=\gamma (a).  \label{4}
\end{equation}

Besides, we have

\begin{equation}
f_{o}=\left\{ 
\begin{array}{c}
\gamma (a-1,1)\text{ \ \ \ for \ \ }a>1 \\ 
\gamma _{o}\text{\ \ \ \ \ \ \ \ \ \ \ \ \ \ \ \ \ \ for \ }a=1%
\end{array}
\right. .  \label{5}
\end{equation}

Hence by (\ref{1}) and (\ref{2}), 
\begin{equation}
\alpha =f_{e}\rho  \label{6}
\end{equation}

\begin{equation}
\alpha ^{\prime }=f_{o}\rho  \label{7}
\end{equation}
and the following relation hold, 
\begin{equation*}
\alpha <f<\alpha ^{\prime }.
\end{equation*}

We shall use the transitive property of adjacency \cite[p. 498]{CLW}, and
negate the possibility of existing two elements $\sigma ,\sigma ^{\prime
}\in L_{n}$, such that $\alpha <\sigma <f$ \ and $f<\sigma ^{\prime }<\alpha
^{\prime }$. Thus, concluding that $\alpha $ is adjacent to $\alpha ^{\prime
}$ in $L_{n}.$

Write $n=md+r$, $0\leq r<m.$

\underline{\textit{Case A}}: For any $r,$ $0\leq r<m,$ if there exist an
element $\sigma \in L_{n}$, such that $\alpha <\sigma <f$ then, we may
choose $\sigma $ in such a way that $\alpha $ is adjacent to $\sigma $ in $%
L_{n}.$ From (\ref{1}) and (\ref{3}), we conclude that 
\begin{equation*}
\sigma =\gamma (a)\rho ^{\prime }
\end{equation*}
since $l(\sigma )>l(f)$.

From (\ref{4}) we have, 
\begin{equation}
\sigma =f_{e}\rho ^{\prime }  \label{7.1}
\end{equation}

Let 
\begin{equation*}
\widetilde{f}=\alpha \wedge \sigma =\gamma (a)(\rho \wedge \rho ^{\prime
})=f_{e}(\rho \wedge \rho ^{\prime }),
\end{equation*}

then, since $f_{e}$ is even, 
\begin{equation}
\widetilde{f}_{e}=f_{e}(\rho \wedge \rho ^{\prime })_{e}  \label{8}
\end{equation}

and 
\begin{equation}
\widetilde{f}_{o}=f_{e}(\rho \wedge \rho ^{\prime })_{o}.  \label{9}
\end{equation}

Suppose that $\widetilde{f}\in L_{\widetilde{m}},$ where obviously, $m<%
\widetilde{m}\leq n$ and let $n=\widetilde{d}\widetilde{m}+\widetilde{r},$ $%
0\leq \widetilde{r}<\widetilde{m}.$

\underline{\textit{Subcase A1}}: $\widetilde{r}>0.$

By Theorem 5.1 \cite[p. 503]{CLW}, since $\alpha $ is adjacent to $\sigma $
in $L_{n}$ then, $\alpha $ and $\sigma $ have the following structure, 
\begin{equation*}
\alpha =(\widetilde{f}_{e})^{\widetilde{d}}\delta
\end{equation*}

and 
\begin{equation*}
\sigma =\widetilde{f}_{o}(\widetilde{f}_{e})^{\widetilde{d}-1}\delta .
\end{equation*}

Substituting from (\ref{6}) and (\ref{8}) we get, 
\begin{equation*}
f_{e}\rho =\{f_{e}(\rho \wedge \rho ^{\prime })_{e}\}^{^{\widetilde{d}%
}}\delta .
\end{equation*}

Equivalently, 
\begin{equation}
\rho =(\rho \wedge \rho ^{\prime })_{e}\{f_{e}(\rho \wedge \rho ^{\prime
})_{e}\}^{^{\widetilde{d}-1}}\delta .  \label{10}
\end{equation}

Similarly, substituting from (\ref{7.1}), (\ref{8})\ and (\ref{9}) we get, 
\begin{equation*}
f_{e}\rho ^{\prime }=f_{e}(\rho \wedge \rho ^{\prime })_{o}\{f_{e}(\rho
\wedge \rho ^{\prime })_{e}\}^{\widetilde{d}-1}\delta .
\end{equation*}

Or equivalently, 
\begin{equation}
\rho ^{\prime }=(\rho \wedge \rho ^{\prime })_{o}\{f_{e}(\rho \wedge \rho
^{\prime })_{e}\}^{\widetilde{d}-1}\delta .  \label{11}
\end{equation}

Letting $\mu =$ $\{f_{e}(\rho \wedge \rho ^{\prime })_{e}\}^{\widetilde{d}%
-1}\delta \ $together with (\ref{4}) and substituting back (\ref{10}) in (%
\ref{6}) and (\ref{11}) in (\ref{7.1}) we get , 
\begin{equation*}
\alpha =\gamma (a)(\rho \wedge \rho ^{\prime })_{e}\mu
\end{equation*}

and 
\begin{equation*}
\sigma =\gamma (a)(\rho \wedge \rho ^{\prime })_{o}\mu .
\end{equation*}

But this says exactly that $\sigma $ is obtained from $\alpha $ by applying
one elementary operation $\mathcal{E}_{s}$ or $\mathcal{E}_{c}$ on the cell
which is the tail of a longer left sequence of even length than $\gamma (a).$
This can happen only if, $\rho $ and $\rho ^{\prime }$ are empty words. But
then, $\alpha =\sigma .$ Hence, there does not exist an element $\sigma \in
L_{n}$, such that $\alpha <\sigma <f$ .

\underline{\textit{Subcase A2}}\textit{:} $\widetilde{r}=0.$

In this case, we must have $\widetilde{d}>1$ and the following relation hold 
\begin{equation*}
\alpha <\widetilde{f}<\sigma <f<\alpha ^{\prime }.
\end{equation*}

By Theorem 4.1, \cite[p. 498]{CLW}, the sequence $\alpha $ has the following
structure, 
\begin{equation*}
\alpha =(\widetilde{f}_{e})^{\widetilde{d}-1}\delta =\widetilde{f}_{e}(%
\widetilde{f}_{e})^{\widetilde{d}-2}\delta .
\end{equation*}

We distinguish between two possibilities.

I. $\widetilde{d}=2.$

Accordingly, $\alpha =$ $\widetilde{f}_{e}\delta $ and by Theorem 4.1, \cite[%
p. 498]{CLW}, class $B_{2}$, we have $\sigma =\widetilde{f}\star \lambda ,$
where $\lambda =(1)$ is the least (and the only) element of $L_{2}.$ Hence, $%
\sigma =\widetilde{f_{o}}\widetilde{f}.$ By Lemma \ref{gosigma}, $\widetilde{%
f}_{o}\delta $ $\in L_{2\widetilde{m}}$. But this says that $\widetilde{f}%
_{o}\delta $ is obtained from $\alpha $ by applying an elementary operation
on the cell which is the tail of a longer left sequence of even length than $%
\gamma (a).$ Hence, there does not exist an element $\sigma \in L_{n}$, such
that $\alpha <\sigma <f$ .

II. $\widetilde{d}\geq 3.$

In this case, one of the following four conditions hold (class $A_{2}$ in
Theorem 4.1, \cite[p. 498]{CLW});

(i) $\widetilde{f}$ is not least in $L_{\widetilde{m}}$, $\widetilde{f_{e}}$
is lexical and $\delta \in L_{\widetilde{m}}$ is adjacent to $\widetilde{f}$
in $L_{\widetilde{m}}.$ By Lemma 2.12, \cite[p. 487]{CLW}, the sequence $%
\widetilde{f_{e}}\delta $ must be lexical because otherwise $\alpha =(%
\widetilde{f}_{e})^{\widetilde{d}-1}\delta ,$ would be nonlexical. By Lemma
2.13, \cite[p.487]{CLW}, and since $\widetilde{f}$ and $\widetilde{f_{e}}%
\delta $ are lexical then, the sequence 
\begin{equation*}
\eta =\widetilde{f}_{o}(\widetilde{f}_{e})^{\widetilde{d}-2}\delta
\end{equation*}

would be lexical where $\eta \in L_{n}$, and 
\begin{equation*}
\alpha <\widetilde{f}<\sigma <\eta <f<\alpha ^{\prime }.
\end{equation*}
Again, this says exactly that $\eta $ is obtained from $\alpha $ by applying
one elementary operation $\mathcal{E}_{s}$ or $\mathcal{E}_{c}$ on the cell
which is the tail of a longer left sequence of even length than $\gamma (a).$
Hence, there does not exist an element $\sigma \in L_{n}$, such that $\alpha
<\sigma <f$ .

(ii) $\widetilde{f}$ is not least in $L_{\widetilde{m}}$, $\widetilde{f_{e}}$
is nonlexical and $\delta \in L_{\widetilde{m}}$ is adjacent to $\widetilde{f%
}$ in $L_{\widetilde{m}}.$ This case is similar to (i).

(iii)$\widetilde{f}$ is not least in $L_{\widetilde{m}}$, $\widetilde{f_{e}}$
is nonlexical and $\delta \in A_{\widetilde{m}-1}$ is nonlexical, but $%
\widetilde{f_{e}}\delta \in L_{2\widetilde{m}}$. Once again, we may apply
the same argument as (i) and conclude that there does not exist an element $%
\sigma \in L_{n}$, such that $\alpha <\sigma <f$ .

(iv) $\widetilde{f}$ is the least sequence in $L_{\widetilde{m}},$ $%
\widetilde{f}=h_{k}(0)\star (2,1^{2(t-1)})$ and $\delta =h_{k}(0)\star
(1^{2t}).$ By Lemma 3.13, \cite[p.495]{CLW}, $\widetilde{f_{e}}\delta \in
L_{2\widetilde{m}}$. Hence, the argument in (i) is applicable and there does
not exist an element $\sigma \in L_{n}$, such that $\alpha <\sigma <f$ .

\underline{\textit{Case B}}: Assume that there exists an element $\sigma
^{\prime }\in L_{n}$, such that $f<\sigma ^{\prime }<\alpha ^{\prime }.$
Then, by subcase A, there is no $\sigma \in L_{n}$ such $\alpha <\sigma <f.$
Thus, we may choose $\sigma ^{\prime }$ such that $\alpha $ is adjacent to $%
\sigma ^{\prime }$ in $L_{n},$%
\begin{equation*}
\alpha <f<\sigma ^{\prime }<\alpha ^{\prime }.
\end{equation*}

From (\ref{2}), (\ref{3}) and (\ref{5}), we must have, 
\begin{equation*}
\sigma ^{\prime }=\left\{ 
\begin{array}{c}
\gamma (a-1,1)\rho ^{\prime },\text{ \ \ for }a>1 \\ 
\gamma _{o}\rho ^{\prime },\text{ \ \ \ \ \ \ \ \ \ \ \ \ \ \ \ for }a=1%
\end{array}
\right. .
\end{equation*}

Or equivalently, 
\begin{equation}
\sigma ^{\prime }=f_{o}\rho ^{\prime }.  \label{12}
\end{equation}

Together with \ (\ref{1}) we get, 
\begin{equation*}
\alpha \wedge \sigma ^{\prime }=f.
\end{equation*}

\underline{\textit{Subcase B1}}: $r>0.$ By Theorem 5.1 \cite[p. 503]{CLW}, $%
\alpha $ and $\sigma ^{\prime }$ have the following structure, 
\begin{equation*}
\alpha =(f_{e})^{d}\delta
\end{equation*}

and 
\begin{equation*}
\sigma ^{\prime }=f_{o}(f_{e})^{d-1}\delta .
\end{equation*}

Substituting, (\ref{6}) and (\ref{12}) we get, 
\begin{equation*}
\rho =(f_{e})^{d-1}\delta
\end{equation*}

and 
\begin{equation*}
\rho ^{\prime }=(f_{e})^{d-1}\delta .
\end{equation*}

Thus, 
\begin{equation*}
\rho =\rho ^{\prime }.
\end{equation*}

But this exactly means from (\ref{7}) and (\ref{12}) that, 
\begin{equation*}
\sigma ^{\prime }=\alpha ^{\prime }.
\end{equation*}

Hence, there is no element $\sigma ^{\prime }\in L_{n}$, such that $f<\sigma
^{\prime }<\alpha ^{\prime }.$

\underline{\textit{Subcase B2}}\textit{: }$r=0.$ By Theorem 5.1 \cite[p. 503]%
{CLW}, we must have $\sigma ^{\prime }$ =$f\star \lambda $ in $L_{n},$ where 
$\lambda $ is the least element of $L_{d}$.
\end{proof}

\bigskip

We may give now a justification of example \ref{example L7}.

\begin{corollary}
If $p$ is prime then, for any $\alpha \in L_{p},$ $\alpha $ is adjacent \ to 
$\alpha ^{\prime }\in L_{p},$ where $\alpha ^{\prime }$ is obtained from $%
\alpha $ by applying one elementary operation $\mathcal{E}_{s}$ or $\mathcal{%
E}_{c}$ on the cell corresponding to the tail of the longest left sequence
of even length, which gives a lexical successor in $L_{p}$.
\end{corollary}

\begin{proof}
Since $p$ is prime then, for any $1<m<p,$ we must have $p=dm+r,$ where $r>0.$
\end{proof}

\bigskip

We state and prove a "dual theorem" of \ref{adjcrth}. The elementary
operations $\mathcal{E}_{s}$ or $\mathcal{E}_{c}$ must be applied on the
cell corresponding to the tail of the longest left sequence of odd length,
which gives a lexical preceding in $L_{n}.$

\begin{theorem}
\label{dualth} Let $\alpha ^{\prime }\in L_{n}$ and suppose that $\alpha
^{\prime }$ is not minimal in $L_{n}.$

\begin{enumerate}
\item If $n=md$ and $\alpha ^{\prime }=g\star \lambda $ for some fundamental
sequence $g\in $ $L_{m}$ then $g_{e}^{m-1}\widehat{g}$ is adjacent to $%
\alpha ^{\prime }$ in $L_{n},$ where $\lambda $ is the least element of $%
L_{d}$ and $\widehat{g}$ is described in the proof of Lemma \ref{gosigma}.

\item Otherwise, $n=md+r,$ $r>0$ and an elementary operation can be applied
on $\alpha ^{\prime }$ to obtain a lexical preceding. Thus, if $\alpha \in $ 
$L_{n}$ is obtained from $\alpha ^{\prime }$ by applying one elementary
operation $\mathcal{E}_{s}$ or $\mathcal{E}_{c}$ on the cell corresponding
to the tail of the longest left sequence of odd length then, $\alpha $ is
adjacent to $\alpha ^{\prime }$ in $L_{n}.$
\end{enumerate}
\end{theorem}

\begin{proof}
By Theorem \ref{adjcrth} we know that each $\alpha ^{\prime }$ occurs
exactly in one of the following two cases;

(1) \ $\alpha ^{\prime }=g\star \lambda $ for some fundamental sequence $%
g\in $ $L_{m}$, where $\lambda $ is the least element of $L_{d}.$ This
happens if and only if $n=md.$ By Theorem 3.4 (2) \cite[p. 1994]{DLW}, $%
g_{e}^{m-1}\widehat{g}$ is adjacent to $\alpha ^{\prime }$ in $L_{n},$ where 
$\widehat{g}$ is described in the proof of Lemma \ref{gosigma}.

(2) Otherwise, $n=md+r,$ $r>0$. Assuming that $\alpha \in $ $L_{n}$ is
obtained from $\alpha ^{\prime }$ by applying one elementary operation on
the cell corresponding to the tail of the longest left sequence of odd
length, then let $\alpha ^{\prime }=f_{o}\delta $ and $\alpha =f_{e}\delta ,$
where $f\in $ $L_{m}.$ If $\alpha $ is not adjacent to $\alpha ^{\prime }$
in $L_{n}$ then there exists a lexical sequence $\beta \in $ $L_{n}$ such
that $\alpha <\beta <\alpha ^{\prime }$and $\beta $ is adjacent to $\alpha
^{\prime }.$ Denote $\widetilde{f}=\beta \wedge \alpha ^{\prime }$ and
suppose that $\widetilde{f}\in L_{\widetilde{m}},$ where $n=\widetilde{m}%
\widetilde{d}+\widetilde{r}.$ It is sufficient to deal with the case $%
\widetilde{r}>0$ since $\widetilde{r}=0$ yields that $\alpha ^{\prime }=%
\widetilde{f}\star \lambda $ and we get back to case (1). Therefore,
assuming $\ \widetilde{r}>0,$ then by Theorem \ref{adjcrth},\ $\beta =%
\widetilde{f_{e}}\sigma $ and $\alpha ^{\prime }=\widetilde{f_{o}}\sigma
=f_{o}\delta $. Since $f_{o}$ is the longest left sequence of \ $\alpha
^{\prime }$ which gives a lexical preceding, we must have $l(\widetilde{f_{o}%
})<l(f_{o}).$ Therefore, $f_{o}=\widetilde{f_{o}}\rho _{e}$ for some non
empty sequence $\rho $ and hence $f=\widetilde{f_{o}}\rho $ while $f_{e}=%
\widetilde{f_{o}}\rho _{o}.$ But then, $\alpha <\beta $ yields $\widetilde{%
f_{o}}\rho _{o}\delta <\widetilde{f_{e}}\sigma $ which gives a contradiction
and the result follows.
\end{proof}

\section{Algorithms for Producing and Ordering the Sets $A_{n},$ $L_{n}$ and 
$D_{n}\label{Algthms}$}

First we emphasize the fact that it is an easy matter to build efficient
algorithm for testing lexicality of any $\alpha \in $ $A_{n}$. This demands
at most $n$ steps needed for comparing $\alpha $ with all its right
sequences.

Besides, one can easily write an efficient algorithm for extracting the
longest left sequence, of even length, which gives a lexical successor. It
needs at most $[\frac{n}{2}]n$ steps. Among them at most $[\frac{n}{2}]$
steps for splitting the negative cells, starting from the last negative
cell, and at most $n$ steps for checking lexicality of the born sequence.

\subsection{The set $A_{n}$ of positive sequences of degree $n$}

There exists an algorithm for producing and ordering lexically all

$\alpha -$sequences in $A_{n}$, starting from any element of $A_{n}.$

\begin{algorithm}
Apply the following steps for producing and ordering lexically all $\alpha -$%
sequences in $A_{n}$:

\begin{enumerate}
\item Start with any element $\alpha \in $ $A_{n}.$

\item Apply one elementary operation $\mathcal{E}_{s}$ or $\mathcal{E}_{c}$
on the tail of the longest left sequence, of $\alpha ,$ of even length and
obtain $\alpha ^{\prime }$, where $\alpha \underset{adj}{<}\alpha ^{\prime
}. $

\item Apply one elementary operation $\mathcal{E}_{s}$ or $\mathcal{E}_{c}$
on the tail of the longest left sequence, of $\alpha ,$of odd length and
obtain $\alpha ^{\prime \prime }$, where $\alpha ^{\prime \prime }\underset{%
adj}{<}\alpha .$

\item Repeat step 2. on $\alpha ^{\prime }$ to obtain all successors of $%
\alpha .$

\item Repeat step 3. on $\alpha ^{\prime \prime }$ to obtain all precedings
of $\alpha .$

\item The algorithm ends once the maximal element $(n)$ and the minimal
element $(1,n-1)$ in $A_{n}$ are obtained.
\end{enumerate}
\end{algorithm}

\begin{proof}
By Lemma \ref{adj-An} and \ since the lexical ordering $"<"$ is total, this
procedure gives out all elements of $A_{n}.$
\end{proof}

\begin{example}
Starting from the sequence $(1,1,1,1)$ we produce the set $A_{4}$ and order
its elements lexically. In each step we have to apply $\mathcal{E}_{s}$ or $%
\mathcal{E}_{c}$ on the cell corresponding to the tail of the longest left
sequence of even length for the elements above $(1,1,1,1),$ and $\ \mathcal{E%
}_{s}$ or $\mathcal{E}_{c}$ on the cell corresponding to the tail of the
longest left sequence of odd length below $(1,1,1,1).$%
\begin{equation*}
\begin{array}{cccc}
+ & - & + & - \\ 
4 &  &  &  \\ 
3 & 1 &  &  \\ 
2 & 1 & 1 &  \\ 
2 & 2 &  &  \\ 
1 & 1 & 2 &  \\ 
1 & 1 & 1 & 1 \\ 
1 & 2 & 1 &  \\ 
1 & 3 &  & 
\end{array}%
\end{equation*}
\end{example}

\subsection{\ \ The set $L_{n}$ of lexical sequences of degree $n$}

There exists an algorithm for producing and ordering all lexical sequences
in $L_{n}.$ We may start from any element as previously, relaying on the
"dual theorem", namely Theorem \ref{dualth}. But Theorem \ref{adjcrth} is
much simpler for application, therefore we start from the least element.

\begin{algorithm}
Apply the following steps for producing and ordering all lexical sequences
in $L_{n}$:

\begin{enumerate}
\item Write $n=2^{l}(2s+1),$ and compute the harmonics of $(0)$ recursively; 
$h_{0}(0)=(0)$, $\ h_{1}(0)=(1)$, $h_{2}(0)=(2,1)$ until $%
h_{l}(0)=[h_{l-1}(0)]_{o}h_{l-1}(0)$

\item Start with the least element $\alpha =\mu \in $ $L_{n}$ where, 
\begin{equation*}
\mu =\left\{ 
\begin{array}{c}
h_{l}(0)\text{ \ \ \ \ \ \ \ \ \ \ \ \ \ \ \ \ \ \ \ \ \ \ \ \ \ \ \ for }s=0
\\ 
h_{l}(0)\star (2,1^{2s-1)})\text{ \ \ \ \ \ \ \ for }s>0%
\end{array}%
\right. .
\end{equation*}

\item Apply one elementary operation $\mathcal{E}_{s}$ or $\mathcal{E}_{c}$
on the tail of the longest left sequence of $\alpha $, of even length, which
gives a lexical successor in $L_{n},$and obtain $\alpha ^{\prime }$.

\item Extract $f=\alpha \wedge \alpha ^{\prime }\in L_{m},$ and divide $\
n=dm+r.$

\item If $r>0$ then, $\alpha \underset{adj}{<}\alpha ^{\prime }$ in $L_{n}.$
Else, if $r=0$ then, write $d=2^{k}(2t+1),$ repeat steps 1.-2. with $l=k$, $%
s=t$ for computing the least element $\lambda $ of $L_{d},$ and obtain $%
\alpha ^{\prime }=f\star \lambda $ which is the adjacent successor, $\alpha 
\underset{adj}{<}f\star \lambda $, in $L_{n}.$

\item Repeat steps 3.-5. on $\alpha ^{\prime }$ to obtain all lexical
successors of $\alpha .$

\item The algorithm ends once the maximal element $(n-1)$ in $L_{n}$ is
obtained.
\end{enumerate}
\end{algorithm}

\begin{proof}
Lemma 3.10 \cite[p. 493]{CLW}, validates the least element formula, whereas,
Theorem \ref{adjcrth} assures that the element $\alpha ^{\prime }$ obtained
from steps 3.-5. is the lexical adjacent successor of $\alpha $. Since the
lexical ordering $"<"$ is total, this procedure gives out all elements in $%
L_{n}.$
\end{proof}

\subsection{The set $D_{n}$}

Recalling that, 
\begin{equation*}
D_{n}=\underset{d\mid n}{\cup }L_{d},
\end{equation*}
we will start from the least element in $L_{n},$ list all its precedings
from $L_{d}$ for $d\mid n$ which are the harmonics of $(0)$ and continue
forward by the steps described in the previous algorithms, taking into
account the harmonics of the fundamental sequence, $f,$ each time we get $%
r=0.$

\begin{algorithm}
Apply the following steps for producing and ordering all lexical sequences
in $D_{n}$:

\begin{enumerate}
\item Write $n=2^{l}(2s+1),$ and compute the harmonics of $(0)$ recursively; 
$h_{0}(0)$, $\ h_{1}(0)$,..., $h_{l}(0).$

\item Compute the least element $\alpha =\mu \in $ $L_{n}$ where, 
\begin{equation*}
\mu =\left\{ 
\begin{array}{c}
h_{l}(0)\text{ \ \ \ \ \ \ \ \ \ \ \ \ \ \ \ \ \ \ \ \ \ \ \ \ \ \ for }s=0
\\ 
h_{l}(0)\star (2,1^{2s-1)})\text{ \ \ \ \ \ \ for }s>0%
\end{array}%
\right. .
\end{equation*}

Then, the first few successive elements in $D_{n}$ are the contiguous
sequences in the following string; 
\begin{equation}
h_{0}(0)=(0)\underset{adj}{<}\ ...\underset{adj}{<}h_{l}(0)\underset{adj}{<}%
\alpha =\mu .  \label{14}
\end{equation}

\item Apply one elementary operation $\mathcal{E}_{s}$ or $\mathcal{E}_{c}$
on the tail of the longest left sequence of $\alpha $, of even length, which
gives a lexical successor in $L_{n},$and obtain $\alpha ^{\prime }$.

\item Extract $f=\alpha \wedge \alpha ^{\prime }\in L_{m},$ and divide $\
n=dm+r.$

\item If $r>0$ then, $\alpha \underset{adj}{<}\alpha ^{\prime }$ in $L_{n}$
\ and extend the string (\ref{14}) in $D_{n}$ as follows: 
\begin{equation*}
h_{0}(0)=(0)\underset{adj}{<}\ ...\underset{adj}{<}h_{l}(0)\underset{adj}{<}%
\alpha =\mu \underset{adj}{<}\alpha ^{\prime }.
\end{equation*}

Else, if $r=0$ then, write $d=2^{k}(2t+1),$ repeat steps 1.-2. with $l=k$, $%
s=t$ for computing the least element $\lambda $ of $L_{d},$ and extend the
string (\ref{14}) in $D_{n}$ as follows: 
\begin{equation*}
h_{0}(0)\underset{adj}{<}\ ...\underset{adj}{<}h_{l}(0)\underset{adj}{<}%
\alpha \underset{adj}{<}f\underset{adj}{<}h_{1}(f)\underset{adj}{<}...%
\underset{adj}{<}h_{k}(f)\underset{adj}{\leq }\alpha ^{\prime }=f\star
\lambda .
\end{equation*}

\item Repeat steps 3.-5. on $\alpha ^{\prime }$ to obtain all lexical
successors of $\alpha .$

\item The algorithm ends once the maximal element $(n-1)$ in $D_{n}$ is
obtained.
\end{enumerate}
\end{algorithm}

\begin{proof}
Here again we use Theorem \ref{adjcrth}, Lemma 3.10 \cite[p. 493]{CLW}, and
Theorem 5.1 \cite[p. 503]{CLW}. Obviously, all elements of $L_{n}$ are
obtained by this algorithm, since the same steps from the previous algorithm
are used. The crucial point is to prove that all elements of $L_{m}$ are
obtained for $m\mid n,$ and $m<n.$ Now, let $\widehat{f}\in L_{m},$ where $%
m\mid n,$ and $m<n.$ If $\widehat{f}$ is not obtained by the algorithm then, 
$\widehat{f}$ must fall between two adjacent elements $\alpha ,$ $\alpha
^{\prime }\in L_{n}.$ Let $f=\alpha \wedge \alpha ^{\prime }$ then clearly
we have, 
\begin{equation*}
\alpha <f<\widehat{f}<\alpha ^{\prime }.
\end{equation*}

By the transitive property of adjacency, $\alpha $ is adjacent to $f$ in $%
L_{n}\cup \left\{ f\right\} $ and $f$ is adjacent to $\alpha ^{\prime }$ in $%
\left\{ f\right\} \cup L_{n}.$

On the other hand, by Theorem 4.1 \cite[p. 498]{CLW}, class B, we have; $%
\widehat{f}$ is adjacent to $\widehat{f}\star \lambda \in L_{n}$ in $\left\{ 
\widehat{f}\right\} \cup L_{n}$ where, $\lambda $ is the least element of $%
L_{d},$ $d=\frac{n}{m}.$

Since the lexical ordering in $D_{n}$ is total then, we have two cases:

\begin{enumerate}
\item[i.] $\alpha <f<\widehat{f}<\widehat{f}\star \lambda <\alpha ^{\prime
}. $

\item[ii.] $\alpha <f<\widehat{f}<\alpha ^{\prime }<\widehat{f}\star \lambda
.$

In case i., $\widehat{f}\star \lambda \in L_{n}$ falls between $f$ and $%
\alpha ^{\prime }$ which is impossible. In case ii., $\alpha ^{\prime }\in
L_{n}$ falls between $\widehat{f}$ and $\widehat{f}\star \lambda $ which is
also impossible. Consequently, all elements of $D_{n}$ are obtained by the
algorithm.
\end{enumerate}
\end{proof}

Finally, we exhibit one more example.

\begin{example}
$D_{8}=L_{1}\cup L_{2}\cup L_{4}\cup L_{8},$ where obviously, $L_{1}=\left\{
(0)\right\} ,$ $L_{2}=\left\{ (1)\right\} ,$ and $L_{4}=\left\{
(2,1),(3)\right\} .$ Now, $8=2^{3}\cdot 1$ so, $l=3$ and $s=0.$ The
harmonics of $(0)$ are $h_{0}(0)=(0),$ $h_{1}(0)=(1)$, $h_{2}(0)=(2,1)$ and $%
h_{3}(0)=(2,1,1,2,1).$ Since $s=0,$ $h_{3}(0)$ is the least element of $%
D_{8}.$ When applying the elementary operations on $(2,1,1,2,1)$ and its
successors, we get only one time $r=0,$ in the case $\alpha =(3,1,2,1),$ $%
\alpha ^{\prime }=(4,1,2).$ Therefore, $f=(3)$ and $h_{1}(f)=(4,3).$ Thus,
applying one more elementary operation on $h_{1}(f)$ we get $\alpha ^{\prime
}=(4,1,2)$ and we have, 
\begin{equation*}
\alpha =(3,1,2,1)\underset{adj}{<}f=(3)\underset{adj}{<}h_{1}(f)=(4,3)%
\underset{adj}{<}\alpha ^{\prime }=(4,1,2).
\end{equation*}%
All elements of $L_{1}\cup L_{2}\cup L_{4}$ are obtained and the full list
of $D_{8}$ is:

\begin{equation*}
\begin{array}{cccccc}
+ & - & + & - & + & - \\ 
7 &  &  &  &  &  \\ 
6 & 1 &  &  &  &  \\ 
5 & 1 & 1 &  &  &  \\ 
5 & 2 &  &  &  &  \\ 
4 & 1 & 2 &  &  &  \\ 
4 & 1 & 1 & 1 &  &  \\ 
4 & 2 & 1 &  &  &  \\ 
4 & 3 &  &  &  &  \\ 
3 &  &  &  &  &  \\ 
3 & 1 & 2 & 1 &  &  \\ 
3 & 1 & 1 & 1 & 1 &  \\ 
3 & 1 & 1 & 2 &  &  \\ 
3 & 2 & 2 &  &  &  \\ 
3 & 2 & 1 & 1 &  &  \\ 
2 & 1 & 2 & 1 & 1 &  \\ 
2 & 1 & 1 & 1 & 1 & 1 \\ 
2 & 1 & 1 & 2 & 1 &  \\ 
2 & 1 &  &  &  &  \\ 
1 &  &  &  &  &  \\ 
0 &  &  &  &  & 
\end{array}%
\end{equation*}
\end{example}

\bigskip

\begin{acknowledgement}
The author is indebted to prof. Arye Juh\`{a}sz from the faculty of
mathematics at the Technion-Israel for numerous discussions and lot of
encouragement for bringing this work out. Also indebted to the anonymous
referee for his valuable comments, in particular for proposing the
guide-line proof of Lemma \ref{gosigma} using reference \cite{DLW}, thus
covering the case $\widetilde{d}=2$ of subcase A2 in Theorem \ref{adjcrth}.
\end{acknowledgement}

\end{document}